\documentclass[10pt, twoside]{article}
\usepackage{amssymb,latexsym,amsmath}
\usepackage[hypertex]{hyperref}
\usepackage[english]{babel}
\begin{document}

\makeatletter
\renewcommand{\@evenfoot}{ \thepage \hfil \footnotesize{\it 
ISSN 1025-6415 \ \ Dopovidi Natsionalno\" \i \ Akademi\" \i \  Nauk Ukra\"\i ni. 2000, no. 6} } 
\renewcommand{\@oddfoot}{\footnotesize{\it 
ISSN 1025-6415 \ \  Dopovidi Natsionalno\" \i \ Akademi\" \i \  Nauk Ukra\"\i ni. 2000, no. 6} 
\hfil \thepage } 
\noindent
\\
\\
\\
\\
\\
\\
\\
\\
\\
\\
\\
\\
\\
\\
\\
\\
\\
\\
\\
\\
\\
\\
\\
\\
\par{
\leftskip=1.6cm  \rightskip=0cm  
\noindent            
UDC 512 \medskip \\
\copyright \ {\bf 2000}\medskip \\
{\bf T. R. Seifullin\medskip \\ }
{\large \bf Koszul complexes of embedded systems of polynomials and 
\smallskip \\
duality \medskip \\ } 
\small {\it (Presented by Corresponding Member of the NAS 
of Ukraine A. A. Letichevsky)} \medskip \\
\small {\it 
The author considers the dependence of Koszul complexes and dependence of dual 
Koszul complexes of two systems of polynomials, when one system is a part 
of other system, in connection with the duality in a Koszul complex 
established by author earlier. Whence, the dependence of Koszul complexes 
and dependence of dual Koszul complexes follow when  one system is linearly
expressed through other system. The obtaned results are used in the proof of
homotopic equiva-lence, formulated earlier by the author, of a Koszul complex 
and a dual Koszul complex, what happens under certain conditions. }
\par \medskip } \noindent 
In the present paper 
%
it is establishing 
the connection between Koszul complexes 
of embedded systems of polynomials under duality, that established  
in the author's paper [3]. 
On the basis of this 
%
it is given 
%
another proof of the connection  Koszul 
complexes of systems of polynomials, connected by linear dependence, 
than has been  earlier obtained by author [4]. The obtained 
result implies theorem 4 in [3].  By novelty of definitions 
and notations, 
for  
their
%
exposition
%
it is require 
a large volume,
for this reason
for all necessary definitions and notations we 
refere the reader 
to works  [3]  and [4],  and also other author's works.
%
%
%
%
%

 Let $a=(a_1,\ldots ,a_n)$  
be a collection of elements, then we will write $|a|=n$.

 Let $({\bf C} ;\partial )$  
be a complex, and 
     $c',c''\in {\bf C} $, 
we will write 
     $c' \buildrel{ \partial }\over \simeq  c''$ 
if 
     $\exists \ c\in {\bf C} $  such  that $c' - c'' = \partial 
      \left[ c\right] $.

 Let ${\bf R}$  
be a commutative ring with unity $1$ and zero $0$, 
if  $x=(x_1,\ldots ,x_n)$  
be a collection of variables, then ${\bf R} [x]$ 
denotes a ring of polynomials in commuting variables $x$ 
with coefficients in ${\bf R} $. 
Let ${\bf A} $  
be a commutative ring with unity $1$ and zero $0$, 
if 
     $\widehat f=(\widehat f_1,\ldots ,\widehat f_s)$  
be a collection of  Grassmann   anticommuting variables, then   
     ${\Lambda} _{\bf A} (\widehat f)$ 
denotes Grassmann algebra in variables $\widehat f$ over ${\bf A}$.

 Let $x=(x_1,\ldots ,x_n)$  and  $y=(y_1,\ldots ,y_m)$  be 
commuting  variables;  
     $\widehat f \ = \ (\widehat f_1,\ldots ,\widehat f_s)$ 
and 
     $\widehat g\ =\ (\widehat g_1,\ldots ,\widehat g_t)$    
be Grassmann anticommuting   variables. Denote by \medskip

\ \ \ ${\bf 1} _{(x,\widehat f)} (a(y),\widehat gb(y))= {\bf 1} _{(x_1,\ldots 
,x_n,\widehat f_1,\ldots ,\widehat f_s)} (a_1(y),\ldots ,a_n(y),\widehat 
gb_1(y),\ldots ,\widehat gb_s(y))=$

\ \ \ $ \qquad = \exp (a(y)x_*+\widehat gb(y)\widehat f_*) = \exp \left( 
\sum\limits_{ k} a_k(y)x^k_*+\sum\limits_{i,j} \widehat g_jb^j_i(y)\widehat 
f^i_*\right)$  
\medskip \\
linear over ${\bf R} $ homomorphism 
     ${\Lambda} _{{\bf R} [x]} (\widehat f) \rightarrow  
     {\Lambda} _{{\bf R} [y]} (\widehat g)$, such that $1 \mapsto 1$,  
     $\forall k=1,n: x_k  \mapsto   a_k(y)$, 
     $\forall i=1,s: \widehat f_i \mapsto  \sum \limits_j\widehat g_jb^j_i(y)$, and we call 
it an exponent.

Let 
     $x=(x_1,\ldots ,x_n)$ and $y=(y_1,\ldots ,y_n)$; $\partial : 
     \widehat u \mapsto  (x-y)  =  (x_1-y_1,\ldots ,x_n-y_n)$; $F(x) 
     \in  {\bf R} [x]$; denote by $\widehat u\nabla F(x,y) = 
     \sum\limits^{ n} _{k=1} \widehat u_k\nabla ^kF(x,y)$ 
such element 
     $\in {\bf C} (x,y,\widehat u)$,  
that \medskip

\ \ \ $\partial \left[ \widehat u\nabla F(x,y)\right] 
   = \partial \left[ \sum\limits^{ n} _{k=1} \widehat u_k\nabla ^kF(x,y)\right] 
 = \sum\limits^{ n} _{k=1} (x_k-y_k)\nabla ^kF(x,y) = F(x) - 
F(y).$ \medskip

 {\it {\bf Lemma 1.}  Let \, ${\bf A}$ \, be \,  a \, commutative \, ring \,  
with \, unity \, $1$ \, and \,  zero \, $0$; 
\ $\widehat f =  (\widehat f_1,\ldots ,\widehat f_s)$,   
$\widehat g = (\widehat g_1,\ldots ,\widehat g_t)$  
be collections of anticommuting variables, 
$\forall k=1,n: \forall i=1,s:  a^i_k\in {\bf A} $  
and
$\forall k=1,n:  \forall j=1,t:$ $b^j_k\in {\bf A} $, then \medskip

\ \ \ $\mathop{\bot}\limits_{ \widehat f} \det \| \widehat f_*\| 
\det \| \widehat fa+\widehat gb\|  = \det \left\| \begin{matrix}a 
& \widehat f_* \cr
\widehat gb &  0\end{matrix} \right\| .$ } \medskip

 {\bf Proof.} \medskip 

\ \ \ $\mathop{\bot}\limits_{ 
\widehat f} \det \| \widehat f_*\| \det \| \widehat fa+\widehat 
gb\|  = \mathop{\top}\limits_{ \widehat f'} \det \| \widehat 
f'_*\| \det \| \widehat f'a+\widehat gb\| \exp (\widehat 
f'\widehat f_*) =$

\ \ \ $ \qquad = \mathop{\top}\limits_{ \widehat f'} \mathop{\top}\limits_{ 
\widehat g'} \exp (\widehat g\widehat g'_*)\det \| \widehat f'_*\| 
\cdot \det \| \widehat f'a+\widehat g'b\| \exp (\widehat f'\widehat 
f_*) = \det \left\| \begin{matrix}a & \widehat f_* \cr
\widehat gb &  0\end{matrix} \right\|. $ \medskip

 {\it {\bf Lemma 2.}  Let \,  ${\bf A} $ \, be \, a \, commutative \,  ring \, 
\, with \, unity \, $1$ \, and \, zero \, $0$; 
$\widehat f  =$ $(\widehat f_1,\ldots ,\widehat f_s)$, 
$\widehat g = (\widehat g_1,\ldots ,\widehat g_t)$  
be collections of anticommuting variables,  
$\forall i=1,s:$ $\forall j=1,t: b^j_i\in {\bf A} $, 
then:

 \ \ \ 1) $\mathop{\bot}\limits_{ \widehat f} \det \| \widehat f_*\| 
\det \| \widehat f-\widehat gb\|  = \det \| \widehat gb\widehat 
f_*\|  = \exp (\widehat gb\widehat f_*)$; \eject

 \ \ \ 2) $\mathop{\top}\limits_{ \widehat f} \det \| \widehat f_*\| 
\det \| \widehat f-\widehat gb\|  = (\widehat g)^0$. }

 {\bf Proof 1.}  By virtue of lemma 1 $$\begin{matrix} 
\mathop{\bot}\limits_{ \widehat f} \det \| \widehat f_*\| 
\det \| \widehat f-\widehat gb\|  & = & \det \left\| \begin{matrix} E_s & \widehat f_* 
\cr
-\widehat gb &  0\end{matrix} \right\|  & = & \det \left\| \begin{matrix} 
E_s &   \widehat 0_* \cr
-\widehat gb & \widehat gb\widehat f_*\end{matrix} \right\|  &  &  \cr 
\vphantom{\left\| \begin{matrix}a \cr
b\end{matrix} \right\| }      &  &        \|    &  &       
  \|     &  &  \cr                     &  & \det \left\| \begin{matrix}E_s 
&   \widehat f_* \cr
\widehat 0 & \widehat gb\widehat f_*\end{matrix} \right\|  & = & \det \left\| 
\begin{matrix}E_s &   \widehat 0_* \cr
\widehat 0 & \widehat gb\widehat f_*\end{matrix} \right\|  & = & \det \| 
\widehat gb\widehat f_*\| ,  \end{matrix}$$ where $\forall i,i'=1,s: 
(E_s)^i_{i'}  = \widehat f^i_*.\widehat f_{i'}  = (i=i')$.

 {\bf Proof 2.}  By virtue of 1) of this lemma \medskip

\ \ \ $\mathop{\top}\limits_{ 
\widehat f} \det \| \widehat f_*\| \det \| \widehat f-\widehat 
gb\|  = \mathop{\top}\limits_{ \widehat f} (\mathop{\bot}\limits_{ 
\widehat f} \det \| \widehat f_*\| \det \| \widehat f-\widehat 
gb\| )(\widehat f)^0 = \mathop{\top}\limits_{ \widehat f} \exp 
(\widehat gb\widehat f_*)(\widehat f)^0 = (\widehat g)^0.$ \medskip 

 {\it {\bf Lemma 3.}  Let \, ${\bf R}$ \, be \, a \, commutative \, ring 
with \, unity \, $1$ \, and \, zero \, $0$; 
     $x=(x_1,\ldots ,x_n)$, $f(x)=(f_1(x),\ldots ,f_s(x)) 
     \in  {\bf R} [x]^s$, $\partial : \widehat f_x \mapsto  
     f(x)$, then \medskip 

\ \ \ $\partial \left[ (x)^0\det \| \widehat f^x_*\| 
\right] = -(f(x)\widehat f^x_*)\det \| \widehat f^x_*\| =0.$ } \\ 

 {\it {\bf Theorem 1.}  Let ${\bf R}$  be a commutative 
ring with unity $1$  and  zero $0$; 
     $x=$ $(x_1,\ldots ,x_n)$, $f(x)=(f_1(x),\ldots ,f_s(x)) 
     \in  {\bf R} [x]^s$, $F(x)=(F_1(x),\ldots ,F_t(x)) 
     \in  {\bf R} [x]^t$; $\partial : \widehat f_x,\widehat f'_x 
     \mapsto  f(x),$ \hbox{$\widehat F_x,\widehat F'_x \mapsto  F(x)$,}  
then \medskip 

\ \ \ $ \partial \left[ x^0\det \| \widehat F'^x_*\| \exp (\widehat 
f_x\widehat f'^x_*)\right]   =  0 \qquad  \& \qquad  \partial 
\left[ x^0(\widehat F_x)^0\exp (\widehat f_x\widehat f'^x_*)\right]   =  
0  \qquad    \qquad   $ \medskip   \\  
and, hence, the maps \medskip

\ \ \ ${\bf C} (x,\widehat f_x,\widehat F_x)   \ni    c(x,\widehat f_x,\widehat F_x) 
  \mapsto    \mathop{\top}\limits_{ \widehat f'_x} \mathop{\top}\limits_{ 
\widehat F'_x} x^0\det \| \widehat F'^x_*\| \exp (\widehat f_x\widehat 
f'^x_*)\, c(x,\widehat f'_x,\widehat F'_x)= $

\ \ \ $\qquad =\mathop{\top}\limits_{ \widehat F_x} \det \| \widehat 
F^x_*\| \, c(x,\widehat f_x,\widehat F_x) \in  {\bf C} (x,\widehat f_x),$

\ \ \ ${\bf C} (x_*,\widehat f^x_*) \ni  c(x_*,\widehat f^x_*) \mapsto  
\mathop{\bot}\limits_{ x} \mathop{\top}\limits_{ \widehat f'_x} x^0\det 
\| \widehat F^x_*\| \exp (\widehat f'_x\widehat f^x_*)\, c(x_*,\widehat 
f'^x_*) =$

\ \ \ $\qquad =\det \| \widehat F^x_*\| \, c(x_*,\widehat f^x_*) 
\in  {\bf C} (x_*,\widehat f^x_*,\widehat F^x_*),$ \medskip

\ \ \ ${\bf C} (x,\widehat f_x) \ni   c(x,\widehat f_x)  \mapsto  
\mathop{\top}\limits_{ \widehat f'_x} x^0(\widehat F_x)^0
\exp (\widehat f_x\widehat f'^x_*)\, c(x,\widehat f'_x) =$

\ \ \ $\qquad = (\widehat F_x)^0 \, c(x,\widehat f_x) \in  {\bf C} (x,\widehat f_x,\widehat 
F_x),$\medskip

\ \ \ ${\bf C} (x_*,\widehat f^x_*,\widehat F^x_*) \ni  
\, c(x_*,\widehat f^x_*,\widehat F^x_*) \mapsto  \mathop{\bot}\limits_{ x} 
\mathop{\top}\limits_{ \widehat f'_x} \mathop{\top}\limits_{ \widehat F'_x} 
x^0(\widehat F'_x)^0\exp (\widehat f'_x\widehat f^x_*)\, 
c(x_*,\widehat f'^x_*,\widehat F'^x_*) =$

\ \ \ $\qquad =\mathop{\top}\limits_{ \widehat F_x} (\widehat F_x)^0
\, c(x_*,\widehat f^x_*,\widehat F^x_*) \in  {\bf C} (x_*,\widehat f^x_*)$ 
\medskip \\
are complex morphisms. } 

 {\bf Proof.}  The equality   
$\partial \left[ x^0\det 
\| \widehat F'^x_*\| \exp (\widehat f_x\widehat f'^x_*)\right] 
=0$   holds, since  \eject\noindent
$\partial \left[ x^0\det \| 
\widehat F'^x_*\| \right] =0$  
(lemma  3) and  
$\partial \left[ x^0\exp (\widehat f_x\widehat f'^x_*)\right] = 0$,
and collections of adjoint variables of factors has empty intersection.

 The equality 
$\partial \left[ x^0(\widehat F_x)^0\exp (\widehat f_x\widehat f'^x_*)\right] =0$  
\ holds, since  
$\partial \left[ x^0(\widehat F_x)^0\right] =0$ 
and 
$\partial \left[ x^0\exp (\widehat f_x\widehat f'^x_*)\right]=0$.\smallskip

 {\it {\bf Theorem 2.}  Let \ ${\bf R}$ \  be a commutative \ ring \ with 
\ unity \ $1$ \ and \ zero \ $0$; 
\ $x=$ $(x_1,\ldots ,x_n)$, $y=(y_1,\ldots ,y_n)$, 
$f(x)=(f_1(x),\ldots ,f_s(x)) 
\in {\bf R} [x]^s$, $F(x)=(F_1(x),\ldots ,F_t(x)) 
\in  {\bf R} [x]^t$; \ $\partial : 
\widehat f_x,\widehat f'_x \mapsto  f(x),  
\ \widehat F_x,\widehat F'_x  \mapsto  F(x)$, 
\ $\widehat f_y,\widehat f'_y  \mapsto $ $f(y), 
\ \widehat F_y,\widehat F'_y \mapsto  F(y), 
\ \widehat u \mapsto  (x-y) = $ \hbox{$(x_1-y_1,\ldots ,x_n-y_n)$}, then 
\medskip

\ \ \ 1) $\mathop{\top}\limits_{ \widehat f_y} \mathop{\top}\limits_{ 
\widehat F_y} \det \left\| \begin{matrix}\nabla F(x,y) & \nabla f(x,y) 
\cr
 \widehat F_x-\widehat F_y &  \widehat f_x-\widehat f_y\end{matrix} \right\| 
y^0\det \| \widehat f^y_*\| \exp (\widehat F'_y\widehat F^y_*)=$
\\
\indent \ \qquad $=\mathop{\top}\limits_{ \widehat F'_x} x^0(\widehat 
f_x)^0\exp (\widehat F_x\widehat F'^x_*)\det \left\| \begin{matrix}\nabla 
F(x,y) \cr
 \widehat F'_x-\widehat F'_y\end{matrix} \right\| $, \medskip \\

\ \ \ 2) $\mathop{\top}\limits_{ \widehat f_y} \det \left\| 
\begin{matrix}\nabla f(x,y) \cr
 \widehat f_x-\widehat f_y\end{matrix} \right\| y^0(\widehat F'_y)^0\exp 
(\widehat f'_y\widehat f^y_*) =$\\
\indent \ \qquad $=\mathop{\top}\limits_{ \widehat f'_x} \mathop{\top}\limits_{ 
\widehat F'_x} (-1)^{| F| | x| } x^0\det \| \widehat F'^x_*\| 
\exp (\widehat f_x\widehat f'^x_*)\det \left\| \begin{matrix}\nabla F(x,y) 
& \nabla f(x,y) \cr
 \widehat F'_x-\widehat F'_y &  \widehat f'_x-\widehat f'_y\end{matrix} \right\| 
$, \bigskip  \\ 
hence, \bigskip

 1') $\mathop{\top}\limits_{ y} \mathop{\top}\limits_{ \widehat 
f_y} \mathop{\top}\limits_{ \widehat F_y} \det \left\| \begin{matrix}\nabla 
F(x,y) & \nabla f(x,y) \cr
 \widehat F_x-\widehat F_y &  \widehat f_x-\widehat f_y\end{matrix} \right\| 
\left( \det \| \widehat f^y_*\| c(y_*,\widehat F^y_*)\right) 
 =$\\
\indent \ \qquad $=(\widehat f_x)^0\left( \mathop{\top}\limits_{ 
y} \mathop{\top}\limits_{ \widehat F_y} \det \left\| \begin{matrix}\nabla 
F(x,y) \cr
 \widehat F_x-\widehat F_y\end{matrix} \right\| c(y_*,\widehat F^y_*)\right) $,
\medskip

 2') $\mathop{\top}\limits_{ y} \mathop{\top}\limits_{ \widehat 
f_y} \det \left\| \begin{matrix}\nabla f(x,y) \cr
 \widehat f_x-\widehat f_y\end{matrix} \right\| \left( \mathop{\top}\limits_{ 
\widehat F_y} (\widehat F_y)^0c(y_*,\widehat f^y_*,\widehat F^y_*)\right)  =$
\\
\indent \ \qquad $=\mathop{\top}\limits_{ \widehat F_x} (-1)^{| F| 
| x| } \det \| \widehat F^x_*\| \left( \mathop{\top}\limits_{ 
y} \mathop{\top}\limits_{ \widehat f_y} \mathop{\top}\limits_{ \widehat 
F_y} \det \left\| \begin{matrix}\nabla F(x,y) & \nabla f(x,y) 
\cr
 \widehat F_x-\widehat F_y &  \widehat f_x-\widehat f_y\end{matrix} \right\| 
c(y_*,\widehat f^y_*,\widehat F^y_*)\right) $. }  \medskip \\ 

 {\bf Proof 1.}  We have \medskip

\ \ \ $ \mathop{\top}\limits_{ 
\widehat f_y} \mathop{\top}\limits_{ \widehat F_y} 
\det \left\| \begin{matrix}\nabla F(x,y) & \nabla f(x,y) \cr
 \widehat F_x-\widehat F_y &  \widehat f_x-\widehat f_y\end{matrix} \right\| 
y^0\det \| \widehat f^y_*\| \exp (\widehat F'_y\widehat F^y_*) 
=$

\ \ \ $\qquad = \mathop{\top}\limits_{ \widehat f_y} 
\det \left\| \begin{matrix}\nabla F(x,y) & \nabla f(x,y) \cr
 \widehat F_x-\widehat F'_y &  \widehat f_x-\widehat f_y\end{matrix} \right\| 
\det \| \widehat f^y_*\|  =$

\ \ \ $\qquad = \mathop{\top}\limits_{ \widehat f_y} \mathop{\top}\limits_{ 
\widehat u} \det \| -\widehat u_*\| 
\det \| \widehat F_x-\widehat 
F'_y-\widehat u\nabla F(x,y)\| \det \| \widehat f_x-\widehat f_y-\widehat 
u\nabla f(x,y)\| \det \| \widehat f^y_*\|  =$

\ \ \ $\qquad = \mathop{\top}\limits_{ \widehat u} \det \| -\widehat 
u_*\| \det \| \widehat F_x-\widehat F'_y-\widehat u\nabla F(x,y)\| 
\left( \mathop{\top}\limits_{ \widehat f_y} \det \| \widehat f_x-\widehat 
f_y-\widehat u\nabla f(x,y)\| \det \| \widehat f^y_*\| \right)  =$

\ \ \ \ $\qquad = \mathop{\top}\limits_{ \widehat u} \det \| -\widehat 
u_*\| \det \| \widehat F_x-\widehat F'_y-\widehat u\nabla F(x,y)\| 
\left( \mathop{\top}\limits_{ \widehat f_y} \det \| \widehat f^y_*\| 
\det \| -\widehat f_x+\widehat f_y+\widehat u\nabla f(x,y)\| \right) =$
\medskip \eject\noindent
(by virtue of 2) of lemma 2)  \smallskip

\ \ \ $\qquad =\mathop{\top}\limits_{ \widehat u} 
\det \| -\widehat u_*\| 
\det \| \widehat F_x-\widehat F'_y-\widehat u\nabla F(x,y)\| \left( (\widehat u)^0(\widehat f_x)^0\right) 
=$

\ \ \ $\qquad =\det \left\| \begin{matrix}\nabla F(x,y) \cr
 \widehat F_x-\widehat F'_y\end{matrix} \right\| (\widehat f_x)^0 = \mathop{\top}\limits_{ 
\widehat F'_x} x^0(\widehat f_x)^0\exp (\widehat F_x\widehat F'^x_*)\det \left\| 
\begin{matrix}\nabla F(x,y) \cr
\widehat F'_x-\widehat F'_y\end{matrix} \right\| $ \medskip

 {\bf Proof 2.}  We have \medskip

\ \ \ $\mathop{\top}\limits_{ 
\widehat f'_x} \mathop{\top}\limits_{ \widehat F'_x} x^0\det \| 
\widehat F'^x_*\| \exp (\widehat f_x\widehat f'^x_*)\det \left\| 
\begin{matrix}\nabla F(x,y) & \nabla f(x,y) \cr
 \widehat F'_x-\widehat F'_y &  \widehat f'_x-\widehat f'_y\end{matrix} 
\right\|$ 
 =

\ \ \ $\qquad = \mathop{\top}\limits_{ \widehat F'_x} \det \| \widehat 
F'^x_*\| \det \left\| \begin{matrix}\nabla F(x,y) & \nabla 
f(x,y) \cr
 \widehat F'_x-\widehat F'_y &  \widehat f_x-\widehat f'_y\end{matrix} \right\| 
 =$

\ \ \ $\qquad = \mathop{\top}\limits_{ \widehat F'_x} \det \| \widehat 
F'^x_*\| \mathop{\top}\limits_{ \widehat u} \det \| -\widehat 
u_*\| \det \| \widehat F'_x-\widehat F'_y-\widehat u\nabla F(x,y)\| 
\det \| \widehat f_x-\widehat f'_y-\widehat u\nabla f(x,y)\|  
=$

\ \ \ $\qquad = (-1)^{| F| | x| } \cdot $

\ \ \ $\qquad \ \cdot \mathop{\top}\limits_{ \widehat u} 
\det \| -\widehat u_*\| ( \mathop{\top}\limits_{ \widehat F'_x} \det 
\| \widehat F'^x_*\| \det \| \widehat F'_x-\widehat F'_y-\widehat 
u\nabla F(x,y)\| ) \det \| \widehat f_x-\widehat f'_y-\widehat 
u\nabla f(x,y)\| =$ \medskip \\ 
(by virtue of 2) of lemma 2) \medskip  

\ \ \ $ \qquad = (-1)^{| F| | x| } \mathop{\top}\limits_{       \widehat 
u} \det \| -\widehat u_*\| \left( (\widehat F'_y)^0(\widehat u)^0\right) 
\det \| \widehat f_x-\widehat f'_y-\widehat u\nabla f(x,y)\|  
=$

\ \ \ $\qquad = (-1)^{| F| | x| } (\widehat F'_y)^0\det \left\| \begin{matrix}\nabla 
f(x,y) \cr
 \widehat f_x-\widehat f'_y\end{matrix} \right\| \, = \, (-1)^{|F| | x|}\,
\mathop{\top}\limits _{ \widehat f_y} \det \left\| \begin{matrix}\nabla 
f(x,y) \cr
 \widehat f_x-\widehat f_y\end{matrix} \right\| y^0(\widehat F'_y)^0\exp 
(\widehat f'_y\widehat f^y_*).$\\  

{\it {\bf Theorem 3.}  Let \ ${\bf R}$ \  be \ a \ commutative \ 
 ring \ with \ unity \ $1$ \ and \ zero \ $0$; 
\ $x= (x_1,\ldots ,x_n)$, \
\ $y= (y_1,\ldots ,y_n)$, \    
$f(x)=(f_1(x),\ldots ,f_s(x)) \in  {\bf R} [x]^s$, \
$F(x)=(F_1(x),\ldots ,F_t(x)) \in  {\bf R} [x]^t$, \
$\forall j=1,t: 
F_j(x) = \sum\limits^{ s} _{i=1} f_i(x)G^i_j(x)$; 
$\partial : \widehat f_x \mapsto  f(x),$ $\widehat F_x 
\mapsto  F(x)$, $\widehat f_y \mapsto  f(y), \widehat F_y \mapsto  F(y)$, 
then \smallskip

\ \ \ $\mathop{\top}\limits_{ 
\widehat F_x} \exp (\widehat f_xG(x)\widehat F^x_*)
\det \left\| \begin{matrix}\nabla F(x,y) \cr
 \widehat F_x-\widehat F_y\end{matrix} \right\|  = \det \left\| \begin{matrix} 
\nabla F(x,y) \cr
\widehat f_xG(x)-\widehat F_y\end{matrix} \right\|  
\buildrel{ \partial}\over \simeq$ 

\ \ \ $\qquad \buildrel{ \partial}\over \simeq
\mathop{\top}\limits_{ \widehat f_y} \det \left\| 
\begin{matrix}\nabla f(x,y) \cr
 \widehat f_x-\widehat f_y\end{matrix} \right\| 
\det \left\| \begin{matrix}G(y) & \widehat f^y_* \cr
-\widehat F_y & 0\end{matrix} \right\| ;$ 
\medskip \\
moreover,  
$\partial 
\left[ \det \left\| \begin{matrix} \nabla F(x,y) \cr
\widehat f_xG(x)-\widehat F_y\end{matrix} \right\| \right]  = 0$, $\partial 
\left[ \det \left\| \begin{matrix}G(y) & \widehat f^y_* \cr
-\widehat F_y & 0\end{matrix} \right\| \right]  = 0.$

 In particular, for 
$L(x_*,\widehat F^x_*) \in  {\bf C} (x_*,\widehat 
F^x_*)$, if $\partial \left[ L(x_*,\widehat F^x_*)\right]  = 0$, 
then \smallskip

\ \ \ $\partial \left[ \mathop{\bot}\limits_{ x} \mathop{\top}\limits_{ 
\widehat F_x} \det \left\| \begin{matrix}G(x) & \widehat f^x_* \cr
-\widehat F_x & 0\end{matrix} \right\| L(x_*,\widehat F^x_*)\right] 
 = 0$ \medskip \\ 
and 
\eject\noindent 
 
\ \ \ $\mathop{\top}\limits_{ \widehat F_x} \exp 
(\widehat f_xG(x)\widehat F^x_*)\mathop{\top}\limits_{ y} 
\mathop{\top}\limits_{ \widehat F_y} 
\det \left\| \begin{matrix}\nabla F(x,y) \cr
 \widehat F_x-\widehat F_y\end{matrix} \right\| L(y_*,\widehat F^y_*) = 
\mathop{\top}\limits_{ y} \mathop{\top}\limits_{ \widehat F_y} \det 
\left\| \begin{matrix} \nabla F(x,y) \cr
\widehat f_xG(x)-\widehat F_y\end{matrix} \right\| L(y_*,\widehat F^y_*) 
\buildrel{ \partial }\over \simeq$

\ \ \ \ $\qquad \buildrel{\partial }\over \simeq  \mathop{\top}\limits_{ 
y} \mathop{\top}\limits_{ \widehat f_y} \det \left\| \begin{matrix}\nabla 
f(x,y) \cr
 \widehat f_x-\widehat f_y\end{matrix} \right\| \left( \mathop{\bot}\limits_{ 
y} \mathop{\top}\limits_{ \widehat F_y} \det \left\| \begin{matrix}G(y) 
& \widehat f^y_* \cr
-\widehat F_y & 0\end{matrix} \right\| L(y_*,\widehat F^y_*)\right).$ } 
\medskip

 {\bf Proof.}  Since 
$\forall j=1,t: \partial 
\left[ -\widehat F_{j,x} +\widehat f_xG_j(x)\right]  = -F_j(x)+f(x)G_j(x) 
=$ $0$, then  

\ \ \ $\partial \left[ \det \| -\widehat F_x+\widehat f_xG(x)\| \right]  = 
\partial \left[ \prod\limits_{j=1,t} (-\widehat F_{j,x} 
+\widehat f_xG_j(x))\right]  = 0$.\medskip

 Let 
$z=(z_1,\ldots ,z_n)$, $\partial : \widehat f_z \mapsto 
 f(z), \widehat F_z \mapsto  F(z)$. 
By virtue of theorem ([3]) \medskip

\ \ \ $\det \left\| \begin{matrix}\nabla F(x,y) & \nabla f(x,y) \cr
 \widehat F_x-\widehat F_y &  \widehat f_x-\widehat f_y\end{matrix} \right\| 
\exp (xz_*+\widehat f_x\widehat f^z_*+\widehat F_x\widehat F^z_*) 
\buildrel{ \partial }\over \simeq $

\ \ \ \ $\qquad \buildrel{ \partial }\over \simeq  \det \left\| 
\begin{matrix}\nabla F(x,y) & \nabla f(x,y) \cr
 \widehat F_x-\widehat F_y &  \widehat f_x-\widehat f_y\end{matrix} \right\| 
\exp (yz_*+\widehat f_y\widehat f^z_*+\widehat F_y\widehat F^z_*). $  
 \medskip  \\
Hence, by  
$\partial \left[ \det \| -\widehat F_z+\widehat 
f_zG(z)\| \right]  = 0$ we have \medskip 

\ \ \ $\det \left\| \begin{matrix}\nabla f(x,y) & \nabla F(x,y) \cr
 \widehat f_x-\widehat f_y &  \widehat F_x-\widehat F_y\end{matrix} \right\| 
\det \| -\widehat F_x+\widehat f_xG(x)\|  =$

\ \ \ \ $\qquad = \mathop{\top}\limits_{ z} \mathop{\top}\limits_{ 
\widehat f_z} \mathop{\top}\limits_{ \widehat F_z} \det \left\| 
\begin{matrix}\nabla F(x,y) & \nabla f(x,y) \cr
 \widehat F_x-\widehat F_y &  \widehat f_x-\widehat f_y\end{matrix} \right\| 
\exp (xz_*+\widehat f_x\widehat f^z_*+\widehat F_x\widehat F^z_*)
\det \| -\widehat F_z+\widehat f_zG(z)\|  \buildrel{ \partial }\over \simeq $

\ \ \ \ $ \qquad \buildrel{ \partial }\over \simeq  \mathop{\top}\limits_{ 
z} \mathop{\top}\limits_{ \widehat f_z} \mathop{\top}\limits_{ \widehat 
F_z} \det \left\| \begin{matrix}\nabla F(x,y) & \nabla f(x,y) 
\cr
 \widehat F_x-\widehat F_y &  \widehat f_x-\widehat f_y\end{matrix} \right\| 
\exp (yz_*+\widehat f_y\widehat f^z_*+\widehat F_y\widehat F^z_*)
\det \| -\widehat F_z+\widehat f_zG(z)\|  =$

\ \ \ \ $ \qquad = \det \left\| \begin{matrix}\nabla F(x,y) & \nabla 
f(x,y) \cr
 \widehat F_x-\widehat F_y &  \widehat f_x-\widehat f_y\end{matrix} \right\| 
\det \| -\widehat F_y+\widehat f_yG(y)\| .$ \medskip \\
Then by 
$\partial \left[ x^0\det \| \widehat F^x_*\| \right] 
 = 0$ and $\partial \left[ y^0\det \| \widehat f^y_*\| \right]  = 0$ 
(lemma 3) 
it holds \medskip

\ \ \ $\mathop{\top}\limits_{ \widehat F_x} \mathop{\top}\limits_{ \widehat f_y} 
x^0\det \| \widehat F^x_*\| \det \left\| 
\begin{matrix}\nabla F(x,y) & \nabla f(x,y) \cr
 \widehat F_x-\widehat F_y &  \widehat f_x-\widehat f_y\end{matrix} \right\| 
y^0\det \| \widehat f^y_*\| \det \| -\widehat F_x+\widehat 
f_xG(x)\|  \buildrel{ \partial }\over \simeq$

\ \ \ \ $\qquad \buildrel{ \partial }\over \simeq  \mathop{\top}\limits_{ 
\widehat F_x} \mathop{\top}\limits_{ \widehat f_y} x^0\det \| \widehat 
F^x_*\| \det \left\| \begin{matrix}\nabla F(x,y) & \nabla f(x,y) \cr
 \widehat F_x-\widehat F_y &  \widehat f_x-\widehat f_y\end{matrix} 
\right\| y^0\det \| \widehat f^y_*\| \det \| -\widehat F_y+\widehat 
f_yG(y)\|.$ 
\medskip \\   
Consider the left part of the equality \medskip

\ \ \ $\mathop{\top}\limits_{ \widehat F_x} \mathop{\top}\limits_{ 
\widehat f_y} x^0\det \| \widehat F^x_*\| \det \left\| \begin{matrix}\nabla 
F(x,y) & \nabla f(x,y) \cr
 \widehat F_x-\widehat F_y &  \widehat f_x-\widehat f_y\end{matrix} \right\| 
y^0\det \| \widehat f^y_*\| \det \| -\widehat F_x+\widehat 
f_xG(x)\|  =$

\ \ \ \ $\qquad =(-1)^{| F| (| F| -| x| )+| F| | F| } \cdot$

\ \ \ \ $ \qquad \cdot \mathop{\top}\limits_{ \widehat F_x} 
\left( \mathop{\bot}\limits_{ \widehat F_x} 
\det \| \widehat F^x_*\| \det \| \widehat F_x-\widehat f_xG(x)\| \right) 
\left( \mathop{\top}\limits_{ \widehat f_y} 
\det \left\| \begin{matrix}\nabla F(x,y) & \nabla f(x,y) \cr 
\widehat F_x-\widehat F_y &  \widehat f_x-\widehat f_y\end{matrix} \right\| 
\det \| \widehat f^y_*\| \right)  = $ \medskip  

\ \ \ (by virtue of 1) of theorem 2 and 1) of lemma 2) 
\eject
\ \ \ \ $\qquad = (-1)^{| F| | x| } 
\mathop{\top}\limits_{ \widehat F_x} \exp (\widehat f_xG(x)
\widehat F^x_*)\left( (\widehat f_x)^0\det \left\| 
\begin{matrix}\nabla F(x,y) \cr
 \widehat F_x-\widehat F_y\end{matrix} \right\| \right)  =$
 
\ \ \ \ $\qquad = (-1)^{| F| | x| } \mathop{\top}\limits_{ \widehat F_x} 
\det \| \widehat f_xG(x)\widehat F^x_*\| \det \left\| \begin{matrix}\nabla 
F(x,y) \cr
 \widehat F_x-\widehat F_y\end{matrix} \right\| , $
\medskip \\
consider the right part of the equality \medskip

\ \ \ $\mathop{\top}\limits_{ \widehat F_x} \mathop{\top}\limits_{ \widehat f_y} 
x^0\det \| \widehat F^x_*\| 
\det \left\| \begin{matrix}\nabla F(x,y) & \nabla f(x,y) \cr
 \widehat F_x-\widehat F_y &  \widehat f_x-\widehat f_y\end{matrix} \right\| 
y^0\det \| \widehat f^y_*\| \det \| -\widehat F_y+\widehat 
f_yG(y)\| =$

\ \ \ \ $ \qquad =\mathop{\top}\limits_{ \widehat f_y} 
\left( \mathop{\top}\limits_{ \widehat F_x} 
\det \| \widehat F^x_*\| \det \left\| \begin{matrix}\nabla 
F(x,y) & \nabla f(x,y) \cr
 \widehat F_x-\widehat F_y &  \widehat f_x-\widehat f_y\end{matrix} \right\| 
\right) \left( \mathop{\bot}\limits_{ \widehat f_y} \det \| 
\widehat f^y_*\| \det \| -\widehat F_y+\widehat f_yG(y)\| \right) = $ \medskip

\ \ (by virtue of 2) of theorem 2 and lemma 1) \medskip 

\ \ \ \ $\qquad =\mathop{\top}\limits_{ \widehat f_y} (-1)^{| F| | x| } \left( 
\det \left\| \begin{matrix}\nabla f(x,y) \cr
 \widehat f_x-\widehat f_y\end{matrix} \right\| (\widehat F_y)^0\right) 
\left( \det \left\| \begin{matrix}G(y) & \widehat f^y_* \cr 
-\widehat F_y & 0\end{matrix} \right\| \right)  =$

\ \ \ \ $\qquad =\mathop{\top}\limits_{ \widehat f_y} (-1)^{| F| | x| 
} \det \left\| \begin{matrix}\nabla f(x,y) \cr
 \widehat f_x-\widehat f_y\end{matrix} \right\| \det \left\| \begin{matrix}G(y) 
& \widehat f^y_* \cr
-\widehat F_y & 0\end{matrix} \right\|.$ \medskip \\ 
Hence, \medskip 

\ \ \ \,$\mathop{\top}\limits_{ \widehat f_y} \det \left\| \begin{matrix}\nabla 
f(x,y) \cr
 \widehat f_x-\widehat f_y\end{matrix} \right\| \det \left\| \begin{matrix}G(y) 
& \widehat f^y_* \cr
-\widehat F_y & 0\end{matrix} \right\|  \buildrel{ \partial }\over 
\simeq  \mathop{\top}\limits_{ \widehat F_x} \exp (\widehat f_xG(x)\widehat 
F^x_*)\det \left\| \begin{matrix}\nabla F(x,y) \cr
 \widehat F_x-\widehat F_y\end{matrix} \right\| .$\bigskip 

 Let us prove theorem 4 in [3]:

 {\it {\bf Theorem 4.}  Let \ ${\bf R}$ \ be \ a \ commutative \ ring \ 
with \ unity \ $1$ \ and \ zero \ $0$; $x=(x_1,\ldots ,x_n)$, 
$f(x)=(f_1(x),\ldots ,f_s(x)) \in  {\bf R} [x]^s$; 
$\partial : \widehat f_x \mapsto  f(x), \widehat f_y \mapsto  f(y)$.  
If the  ring   
${\bf H} _0(x,\widehat f_x) \simeq {\bf R} [x]/(f(x))_x$   
is a finitely generated module over ${\bf R} $, then } \medskip

\ \ \ $\exists \ e(x_*,\widehat f^x_*) 
\in  {\bf Z} (x_*,\widehat f^x_*): \mathop{\top}\limits_{ y} 
\mathop{\top}\limits_{ \widehat f_y} 
\det \left\| \begin{matrix}\nabla f(x,y) \cr 
\widehat f_x-\widehat f_y\end{matrix} \right\| e(y_*,\widehat f^y_*) 
\buildrel{ \partial }\over \simeq  x^0(\widehat f_x)^0 = 1.$ \medskip

 {\bf Proof.} From commutative algebra it is known, 
that if the ring ${\bf R} [x]/(f(x))_x$ is a finitely generated module over 
${\bf R} $,  then  for  any polynomial $h(x)\in {\bf R} [x]$ 
there exists a monic, i. e. 
with coefficient $1$ of the greatest degree of variable, 
polynomial 
$T(h) = (h)^d+a_{d-1} (h)^{d-1} +\ldots +a_0(h)^0 \in {\bf R} [h]$  
such that  $T(h(x)) \in (f(x))_x$.  

This is proved as follows:
we take the characteristic polynomial of the matrix 
of the multiplication operator by 
$h(x)$ on  ${\bf R} [x]/(f(x))_x$ 
in some finite system of generators of ${\bf R} [x]/(f(x))_x$ as 
a module over  ${\bf R} $,  
this is a monic polynomial, a root of which is $h(x)$.

 Then $\forall j=1,n:$ there exist a monic polynomial 
$T_j(x_j)\in {\bf R} [x_j]$ such  that $T_j(x_j)$ $\in  (f(x))_x$, 
hence,
$T_j(x_j) = \sum\limits^{ s} _{i=1} f_i(x)G^i_j(x)$. Let $\partial 
: \widehat T_{j,x_j}   \mapsto   T_j(x_j)$, $\widehat T_{j,y_j}  \mapsto 
 T_j(y_j)$, 
then there exists 
$L_j(x^j_*,\widehat T^{j,x_j} _*) = 
l_j(x^j_*){\bf 1} _{\widehat T_{j,x_j} } (0) \in  {\bf C} (x^j_*,\widehat 
T^{j,x_j} _*)$ 
such that \medskip   

\ \ \ $\partial \left[ L_j(x^j_*,\widehat 
T^{j,x_j} _*)\right]     =    \mathop{\bot}\limits_{ x_j} T_j(x_j)l_j(x^j_*)\widehat 
T^{j,x_j} _*    =    0$    
\eject\noindent
and \medskip 

\ \,$\mathop{\top}\limits_{ y_j} \mathop{\top}\limits_{ 
\widehat T_{y_j} } \det \left\| \begin{matrix} \nabla T_j(x_j,y_j) 
\cr
\widehat T_{j,x_j} -\widehat T_{j,y_j} \end{matrix} \right\| 
L_j(y^j_*,\widehat T^{j,y_j} _*)
=\mathop{\top}\limits_{ y_j} \det \left\| \nabla T_j(x_j,y_j)\right\| 
l_j(y^j_*)(\widehat T_{j,x_j} )^0= (x_j)^0(\widehat T_{j,x_j} )^0.$
\medskip \\
Indeed, the condition 
%
$\mathop{\bot}\limits_{ x_j} T_j(x_j)l_j(x^j_*) = 0$ 
is equivalent to the following condition: \medskip

\ \ $\forall {\delta} \geq 0: 
\mathop{\top}\limits_{ x_j} (x_j)^{\delta} T_j(x_j)l_j(x^j_*) =$

\ \ \ \ $\qquad = \mathop{\top}\limits_{ x_j} (x_j)^{\delta} \left( 
(x_j)^{d_j} +a_{j,d_j-1} (x_j)^{d_j-1} +\ldots +a_{j,0} (x_j)^0\right) 
l_j(x^j_*) =$

\ \ \ \ $\qquad = l_j(x^j_*).(x_j)^{d_j+{\delta} } +a_{j,d_j-1} 
( l_j(x^j_*).(x_j)^{d_j+{\delta} -1} ) +\ldots +a_{j,0} 
( l_j(x^j_*).(x_j)^{\delta} )   = 0, $
\medskip \\ 
then any  values of 
$({\lambda} _0,\ldots ,{\lambda} _{d_j-1} )\in {\bf R} ^{d_j} $  
uniquely determine $l_j(x^j_*)$ such that 
$\forall {\delta} =0,d_j-1: l_j(x^j_*).(x_j)^{\delta} 
 = {\lambda} _{\delta} $ and $\mathop{\bot}\limits_{ x_j} T_j(x_j)l_j(x^j_*) 
 = 0$. If $l_j(x^j_*)$ such that $\forall {\delta} =0,d_j-2: 
l_j(x^j_*).(x_j)^{\delta}  = 0$ and $l_j(x^j_*).(x_j)^{d_j-1}  
 =  1$, then $\mathop{\top}\limits_{ y_j} \det \left\| \nabla T_j(x_j,y_j)\right\| 
l_j(y^j_*) = (x_j)^0$.

 Let $F(x)=(F_1(x),\ldots ,F_n(x))$,  where   
$\forall j=1,n: F_j(x) = T_j(x_j) =$ 
$\sum\limits^{ s} _{i=1} f_i(x)G^i_j(x)$, 
and  let  
$\partial :  \widehat F_x  \mapsto   F(x),  \widehat F_y  \mapsto 
  F(y)$,  then  $L(x_*,\widehat F^x_*)  =$ $\prod\limits_{j=1,n} 
L_j(x^j_*,\widehat F^{j,x} _*) = l(x_*){\bf 1} _{\widehat F_x} (0)$ such 
 that $\partial \left[ L(x_*,\widehat F^x_*)\right]  = 0$,  since 
 $\forall j=1,n:$ $\partial \left[ L_j(x^j_*,\widehat F^{j,x} _*)\right] = 0$,
and \medskip 

\ \ \ $\mathop{\top}\limits_{ y} \mathop{\top}\limits_{ 
\widehat F_y} \det \left\| \begin{matrix}\nabla F(x,y) \cr
 \widehat F_x-\widehat F_y\end{matrix} \right\| L(y_*,\widehat F^y_*)  
 = \left( \mathop{\top}\limits_{ y} \det \left\| \nabla F(x,y)\right\| 
l(y_*)\right) (\widehat F_x)^0 = (x)^0(\widehat F_x)^0, $ 
\medskip \\
since \medskip

\ \ \ $\mathop{\top}\limits_{ y} \det \left\| \nabla 
F(x,y)\right\| l(y_*) = 
\mathop{\top}\limits_{ y_1} \ldots \mathop{\top}\limits_{ y_n}  
\left( \prod\limits_{j=1,n} \det \| \nabla T_j(x_j,y_j)\| 
\right) \left( \prod\limits_{j=1,n} l_j(y^j_*)\right)=$

\ \ \ $\qquad = \prod\limits_{j=1,n} \left( \mathop{\top}\limits_{ 
y_j} \det \left\| \nabla T_j(x_j,y_j)\right\| l_j(y^j_*)\right) 
 = \prod\limits_{j=1,n} (x_j)^0 = (x)^0.$   
\medskip \\
Denote by \medskip

\ \ \ $e(x_*,\widehat f^x_*) = \mathop{\bot}\limits_{ x} 
\mathop{\top}\limits_{ \widehat F_x} \det \left\| \begin{matrix}G(x) 
& \widehat f^x_* \cr
-\widehat F_x & 0\end{matrix} \right\| L(x_*,\widehat F^x_*) =$

\ \ \ $\qquad =\mathop{\bot}\limits_{ x} \mathop{\top}\limits_{ \widehat F_x} 
\det \left\| \begin{matrix}G(x) & \widehat f^x_* \cr
-\widehat F_x & 0\end{matrix} \right\| l(x_*){\bf 1} _{\widehat F_x} 
(0) = \mathop{\bot}\limits_{ x} \det \| \begin{matrix}G(x) 
& \widehat f^x_*\end{matrix} \| l(x_*),$  
\medskip \\ 
then by virtue of theorem 3 
it hold: 
$\partial \left[ e(x_*,\widehat f^x_*)\right] = 0$ and \medskip

\ \ \ $\mathop{\top}\limits_{ y} 
\mathop{\top}\limits_{ \widehat f_y} 
\det \left\| \begin{matrix}\nabla f(x,y) \cr 
\widehat f_x-\widehat f_y\end{matrix} \right\| 
e(y_*,\widehat f^y_*)= \mathop{\top}\limits_{ y} 
\mathop{\top}\limits_{ \widehat f_y} 
\det \left\| \begin{matrix}\nabla f(x,y) \cr 
\widehat f_x-\widehat f_y\end{matrix} 
\right\| \left( \mathop{\bot}\limits_{ y} 
\mathop{\top}\limits_{ \widehat F_y} 
\det \left\| \begin{matrix}G(y) & \widehat f^y_* \cr
-\widehat F_y & 0\end{matrix} \right\| 
L(y_*,\widehat F^y_*)\right) \buildrel{ \partial }\over \simeq$  
\eject

\ \ \ $\qquad  
\buildrel{ \partial }\over \simeq  \mathop{\top}\limits_{ \widehat 
F_x} \exp (\widehat f_xG(x)\widehat F^x_*)\mathop{\top}\limits_{ y} \mathop{\top}\limits_{ 
\widehat F_y} \det \left\| \begin{matrix}\nabla F(x,y) \cr
 \widehat F_x-\widehat F_y\end{matrix} \right\| L(y_*,\widehat F^y_*) =$ 

\ \ \ $\qquad = \mathop{\top}\limits_{ \widehat F_x} \exp (\widehat f_xG(x)\widehat 
F^x_*)(x)^0(\widehat F_x)^0 = (x)^0(\widehat f_x)^0.$ 
\medskip \\

{\footnotesize

\begin{enumerate}

\item {\it Seifullin, T. R.} 
Root functionals and root polynomials 
of a system of polynomials. (Russian)
Dopov. Nats. Akad. Nauk Ukra\"\i ni  -- 1995, -- no. 5, 5--8.

\item {\it Seifullin, T. R.} Root functionals and root relations 
of a system of polynomials. (Russian) 
Dopov. Nats. Akad. Nauk Ukra\"\i ni  -- 1995, -- no. 6, 7--10.

\item  {\it Seifullin, T. R.}  Homology of the Koszul complex of a 
system of polynomial equations. (Russian)
Dopov. Nats. Akad. Nauk Ukr. Mat. Prirodozn. Tekh. Nauki 1997, no. 9, 43--49. 
\href{http://arxiv.org/abs/1205.0472} {{\tt  arXiv:1205.0472}} (English).

\item  {\it Seifullin, T. R.}  Koszul complexes of systems of 
polynomials connected by linear dependence. (Russian) 
Some problems in 
contemporary mathematics (Russian), 326--349, Pr. Inst. Mat. Nats. Akad. Nauk 
Ukr. Mat. Zastos., 25, Natsional. Akad. Nauk Ukra\"\i ni, Inst. Mat., Kiev, 
1998.
\\
\end{enumerate}

\small{\noindent
{\it V. M. Glushkov Institute of Cybernetics of the NAS of Ukraine, Kiev
\hfill Received 23.02.99
\medskip  \\ 
E-mail: \ {\tt  timur\_sf@mail.ru}
\\
}
}

\end{document}